\documentclass[11pt]{article}
\usepackage{amssymb}

\newtheorem{thm}     {Theorem}[section]
\newtheorem{prop}    [thm]{Proposition}

\newtheorem{lemma}   [thm]{Lemma}

\newcommand{\proof} {\noindent{\bf Proof. }}

\newcommand{\B}{\mathbb B}
\newcommand{\C}{\mathbb C}
\newcommand{\D}{\mathbb D}

\newcommand{\R}{\mathbb R}

\newcommand{\PP}{\mathbb P}

\begin{document}

\title{ Holomorphic mappings between domains with low boundary regularity}
\author{Alexandre Sukhov }
\date{}
\maketitle

\bigskip
 
{\small Abstract.  We study the boundary regularity of proper holomorphic mappings between strictly pseudoconvex domains with boundaries of low regularity.}

MSC: 32H02.

Key words: strictly pseudoconvex domain, proper holomorphic mapping, boundary regularity.

\bigskip

{\small

 Universit\'e  de Lille, Laboratoire
Paul Painlev\'e,
U.F.R. de
Math\'ematiques, 59655 Villeneuve d'Ascq, Cedex, France, e-mail: sukhov@math.univ-lille1.fr,
 and  Institut of Mathematics with Computing Centre , Ufa Federal Research Centre of Russian
Academy of Sciences, 450077, Chernyshevsky Str. 112, Ufa, Russia.

The author is partially suported by Labex CEMPI.
}

\section{Introduction}

The present paper considers an old problem of precise boundary regularity of a proper holomorphic mapping between two strictly pseudoconvex domains in the case when at least one of the boundaries is of regularity  exactly $C^2$. We prove the following

\begin{thm}
\label{MainTheo1}
Let $\Omega_1$ and $\Omega_2$ be bounded strictly pseudoconvex domains in $\C^n$. Suppose that the boundary of $\Omega_1$ is of class $C^{2+ \varepsilon}$ with $\varepsilon > 0$ and the boundary of $\Omega_2$ is of class $C^2$.   Assume that 
 $f:\Omega_1 \to \Omega_2$ is a proper holomorphic mapping. Then $f$ extends to a mapping of class $C^{\alpha}(\overline{\Omega}_1)$
for each $\alpha \in [0,1[$.
\end{thm}

At present the boundary regularity of proper or biholomorphic mappings between strictly pseudoconvex domains is very well understood. Ch. Fefferman \cite{F} proved that a biholomorphic mapping betweeen strictly pseudoconvex domains with boundaries of class $C^\infty$ extends as a $C^\infty$ diffeomorphism between their closures. His proof is based on the study of asymptotic behavior of the Bergman kernel near the boundary. Later several different approaches have been explored. Some of them allow to study the case where the boundaries of domains are of finite smoothness. S.Pinchuk and S.Khasanov \cite{PiKh} proved that a proper holomorphic mapping between strictly pseudoconvex domains with $C^s$ boundaries, with $s > 2$  real, extends to the boundary as a mapping of class $C^{s-1}$ if $s$ is not an integer, and of class $C^{s-1-\varepsilon}$, with any $\varepsilon > 0$, when $s$ is an integer. Y. Khurumov proved that a similar result still remains true   with the loss of regularity on $1/2$. 
The question on the precise regularity in the case where at least one of domains has the boundary of class $C^2$, remains open for a long period. Y.Khurumov announced without any details that his result remains true also in this case, but, to the best of my knowledge, a detailed  proof is not available. The only well-known result (see, for instance, \cite{Pi1}) states that a mapping extends to the boundary as a H\"older $1/2$-continuous mapping. Theorem \ref{MainTheo1} makes a first step toward a definitive answer.
Our main tool is the result of E.Chirka, B. Coupet and the author \cite{ChSu} giving a precise boundary regularity (in the H\"older scale) of a complex disc with boundary glued to a totally real manifold of class $C^1$.  This allows to improve the boundary regularity of a mappping. Note that Theorem \ref{MainTheo1} is new even in the case where the boundary of $\Omega_1$ is of class $C^\infty$. It is quite expectable that the method of the present paper also allows to deal with the case where the boundary of $\Omega_1$ is exactly of class $C^2$. In the last section we sketch suitable modifications of our agument.

The present paper is devoted to the memory of A.G.Vitushkin. It was written  when the author visited the Mathematics Department of the Indiana University (Bloomington) for the Spring semester of 2020. I thank this institution for excellent conditions of the work.

\section{Preliminaries}

We briefly recall some well-known definitions and basic notations.

\subsection{Classes of domains and functions} Let $\Omega$ be a domain in $\C^n$. For a positive integer $k$, denote by  $C^k(\Omega)$  the space of $C^k$-smooth complex-valued functions in~$\Omega$. Also $C^k(\overline\Omega)$ denotes the class of functions whose partial derivatives up to order $k$ extend as continuous functions on 
$\overline\Omega$. Let  $s > 0$ be a real noninteger and $k$ is its integer part. Then  $C^s(\Omega)$ denotes the space of functions of class $C^k(\overline\Omega)$ such that their partial derivatives of order $k$ are (global)  $(s-k)$-H\"older-continuous in $\Omega$; these derivatives  automatically satisfy the  H\"older condition on $\overline\Omega$ so the notation $C^s(\overline\Omega)$ for the same space of functions is also appropriate.


A (closed) real submanifold $E$ of a domain $\Omega \subset \C^n$ is of class $C^s$ (with real $s  \ge 1$) if for every point $p \in E$ there exists an open neighbourhood $U$ of $p$ and a map $\rho: U \longrightarrow \R^d$ of the maximal rank $d<2n$ and of class $C^s$  such that $E \cap U = \rho^{-1}(0)$; then $\rho$ is called a local defining (vector-valued ) function of $E$. The positive integer $d$ is the real codimension of $E$. In the most important special case $d=1$ we obtain the class of real hypersurfaces.

Let $J$ denotes the standard complex structure of $\C^n$. In other words, $J$ acts on a vector $v$ by multiplication by $i$ that is $J v = i v$. For every $p \in E$ the {\it holomorphic tangent space}  $H_pE:= T_pE \cap J(T_pE)$ is the maximal complex subspace of the tangent space $T_pE$ of $E$ at $p$. Clearly $H_pE =   \{ v \in \C^n: \partial \rho(p) v = 0 \}$. The complex dimension of $H_pE$ is called the CR dimension of $E$ at $p$; a manifold $E$ is called a {\it CR  (Cauchy-Riemann) manifold} if its CR dimension is independent of $p \in E$.

A real submanifold $E \subset \Omega$ is called {\it generic} (or generating) if the complex span of $T_pE$ coincides with $\C^n$ for all $p \in E$. Note that every generic manifold of real codimension $d$ is a CR manifold of CR dimension $n-d$. A function $\rho = (\rho_1,...,\rho_d)$ defines a generic manifold if $\partial\rho_1 \wedge ...\wedge \rho_d \neq 0$.
Of special importance are the so-called {\it totally real manifolds}, i.e., submanifolds $E$ for which $H_pE = \{ 0 \}$ at every $p \in E$. A totally real manifold in $\C^n$  is generic if and only if its real dimension is equal to $n$; this is the maximal possible value for the dimension of a totally real manifold. 

 Let $\Omega$ be a bounded domain in~$\C^n$. Suppose that its boundary $b\Omega$ is a (compact) real hypersurface of class $C^s$ in $\C^n$. Then there exists a $C^s$-smooth real 
function $\rho$ in a neighbourhood $U$ of the closure $\overline\Omega$ such that $\Omega = \{ \rho < 0 \}$ and 
$d\rho|_{b\Omega} \ne 0$. We call such a function $\rho$ a global defining function. If $s \ge 2$ one may consider 
{\it the Levi form} of $\rho$:
\begin{eqnarray}
\label{Leviform}L(\rho,p,v) = \sum_{j,k = 1}^n \frac{\partial^2\rho}{\partial z_j\partial\overline{z}_k}(p)v_j \overline v_k.
\end{eqnarray}
A bounded domain $\Omega$ with $C^2$ boundary is called  {\it strictly pseudoconvex}
if $L(\rho,p,v) > 0$  for every nonzero  vector $v \in H_p(b\Omega)$.

\subsection{The Kobayashi-Royden pseudometric}  Denote by $\D = \{ \zeta \in \C: \vert \zeta \vert < 1 \}$ the unit disc in $\C$. We also denote by $\B$ the unit ball of $\C^n$ (the dimension $n$ will be clear from context).  Let $\Omega$ be a domain in $\C^n$; denote by $\mathcal O (\mathbb D, \Omega)$ the class of holomorphic maps from $\D$ to $\Omega$.

Let $z$ be a point of a domain $\Omega$ and $v$ be a tangent vector at $z$. The infinitesimal Kobayashi-Royden pseudometric $F_\Omega(z,v)$  is defined as 
\begin{equation}
F_\Omega(z,v) = \inf \left\{\lambda>0: \exists\, h \in \mathcal O (\mathbb D, \Omega) {\rm \ with \ } 
h(0)=z,\ h'(0)=\frac{v}{\lambda} \right\} .
\end{equation}
This is a nonnegative upper semicontinuous function on the tangent bundle of $\Omega$; its integrated form coincides with the usual Kobayashi distance. The Kobayashi-Royden metric is decreasing under holomorphic mappings: if $f: \Omega \to \Omega'$ is a holomorphic mapping between two domains in $\C^n$ and $\C^m$ respectively, then
\begin{eqnarray}
\label{Kob1}
F_{\Omega'}(f(z),df(z)v) \le F_\Omega(z,v) .
\end{eqnarray}
In fact, this is the largest metric in the class of  infinitesimal metrics that are decreasing under holomorphic mappings.
It is easy to obtain an upper bound on $F_\Omega$. Indeed, let $z + R\B$ with $R = {\rm dist\,}(z,b\Omega)$ be the ball
contained in $\Omega$. It follows by the holomorphic decreasing property applied to the natural inclusion 
$\iota: z + R\B  \to \Omega$ that the Kobayashi-Royden metric of this ball is bigger than $F_\Omega$. This gives the upper bound
\begin{eqnarray}
\label{Kob2}
F_\Omega(z,v) \le \frac{C \vert v \vert}{{\rm dist\,} (z,b\Omega)} .
\end{eqnarray}
Lower bounds require considerably more subtle analysis. They are obtained by several authors employing  various methods.  Quite  general estimates can be established using plurisubharmonic functions. 



For example, this approach leads to the following result inspired by the work of  N.Sibony \cite{Si}:

\begin{prop}
\label{KobProp1}
Let $\Omega$ be a domain in $\C^n$ and $\rho$ be a negative $C^2$-smooth plurisubharmonic function in $\Omega$. Suppose that the partial derivatives of $\rho$ are bounded on $\Omega$ and there exists a constant $C_1 > 0$ such that 
\begin{eqnarray}
\label{Kob4_1}
L(\rho,z,v) \ge C_1 \vert v \vert^2
\end{eqnarray}
for all $z$ and $v$. Then there exists a constant $C_2 > 0$, depending only on the $C^2$-norm of $\rho$, such that 
\begin{eqnarray}
\label{Kob4}
F_\Omega(z,v) \ge C_2 \left ( C_1^2\frac{\vert \langle \partial \rho(z), v \rangle \vert^2}{\vert \rho(z) \vert^2} + C_1\frac{\vert v \vert^2}{\vert \rho(z) \vert^2} \right ).
\end{eqnarray}
\end{prop}
Note that $\rho$ is not assumed to be a defining function of $\Omega$, although this special case is particularly important in applications. The original argument of Sibony assumes that $\Omega$ is globally bounded but this condition can be dropped. In fact, the estimate (\ref{Kob4}) holds on an open subset of $\Omega$ where (\ref{Kob4_1}) is satisfied. Therefore, it can be used in order to localize the Kobayashi-Royden metric.  Note also that $\Omega$ is not assumed to be bounded or hyperbolic. Of course, for the present paper the bounded case is sufficient.

The following localization principle for the Kobayashi-Royden pseudometric is established in \cite{ChSu}. The proof also is based on methods of the prluripotential theory.

\begin{prop}
\label{ProLoc1}
Let $D$ be a domain in $\C^n$. Suppose that $u$ is a negative plurisubharmonic function in $D$ such that the 
for some constants $\varepsilon, B > 0$ the following holds:
\begin{itemize}
\item[(i)] the function $u(z) - \varepsilon \vert z \vert^2$ is plurisubharmonic on 
$D \cap 3 \B$;
\item[(ii)] $\vert u \vert \le B$ on $2 \B$.
\end{itemize}
Then there exists a positive constant $M = M(\varepsilon,B)$, independent of $u$, such that 
$$F_D(w,\xi) \ge M \vert \xi \vert \vert u(w) \vert^{-1/2}$$
for $w \in D \cap 2 \B$.
\end{prop}
Note that this result  also gives a lower bound for the metric. This will be  considerably used in our proof.

\subsection{Complex discs} Consider a wedge-type domain 

\begin{eqnarray}
\label{wedge}
W = \{ z \in \C^m: \phi_j(z) < 0, j = 1,...,m \}
\end{eqnarray}
with the edge 

\begin{eqnarray}
\label{edge}
E = \{ z \in \C^m: \phi_j(z) = 0, j = 1,...,m \}
\end{eqnarray}
We assume that the defining functions $\phi_j$ are of class $C^{1 + \alpha}$ with $\alpha > 0$. Furthermore, as usual we suppose that $E$ is a generic manifold that is 
$\partial \phi_1 \wedge ... \wedge \partial \phi_m \neq 0$ in a neighborhood of $E$. 

Given $\delta > 0$ (which is supposed to be small enough) we also define a shrinked wedge 
\begin{eqnarray}
\label{shrinked}
W_{\delta} = \{ z \in \C^m: \phi_j - \delta \sum_{l \neq j} \phi_l < 0, j= 1,...,m \} \subset W
\end{eqnarray}

We need the well-known construction of filling a wedge $W$ by complex discs gluing to $E$ along an open arc.  A complex disc is a holomorphic map $h: \D \to \C^n$ which is at least contious on the closed disc $\overline\D$. Denote by $b\D^+$ the upper semi-circle. 

\begin{prop}
\label{discs}
Fix $\delta > 0$. There exists a map $H: \D \times \R^{m-1} \to W$ of class $C^{1 + \alpha}(\overline \D \times \R^{m-1})$, $H:( \zeta,t) \mapsto h_t(\zeta)$  with the following  properties:
\begin{itemize}
\item[(i)] for every $t \in \R^{m-1}$ the map $\zeta \mapsto h_t(\zeta)$ is holomorphic in $\D$ and $h(b\D^+)$ is contained in $E$.
\item[(ii)] the curves $h_t(b\D^+)$ from a foliation of $E$.
\item[(iii)]  every disc $h_t(\overline\D)$ is transverse to $E$.
\item[(iv)] $W_\delta \subset \cup_t h_t(\D)$.
\end{itemize}
\end{prop}

The above  {\it gluing disc argument} is often quite helpful for the study of totally real submanifolds. We sketch the idea of proof.  It was introduced in \cite{Pi1} and then used by many authors. Without loss of generality, we may assume that in a neighbourhood $\Omega$ of the origin a smooth totally real manifold $E$ is defined by the equation $x = r(x,y)$, where a smooth vector function $r = (r_1,...,r_n)$ satisfies the conditions $r_j(0) = 0$, $d r_j(0) = 0$. Fix a positive noninteger $s$ and consider for a real function $u \in C^s(b\D)$  the Hilbert transform $T: u \to T(u)$. It is uniquely defined by the conditions that  the function $u + iT(u)$ is the trace of a function holomorphic on $\D$ and $T(u)$ vanished at the origin. Explicitely it is given  by the singular integral 

\begin{eqnarray*}
T(u)(e^{i\theta}) = \frac{1}{2\pi} v. p. \int_{-\pi}^{\pi} u(e^{it}) \cot \left ( \frac{\theta - t}{2} \right ) dt
\end{eqnarray*}

 This is a classical linear singular integral operator; it is  bounded on the space $C^s(b\D)$. Let $S^+ = \{ e^{i\theta}: \theta \in [0,\pi] \}$ and $S^- = \{ e^{i\theta}: \theta \in ]\pi, 2 \pi[ \}$ be the semicircles. Fix a $C^\infty$-smooth real function $\psi_j$ on $b\D$ such that $\psi_j\vert S^+ = 0$ and  $\psi_j\vert S^- < 0$, $j=1,...,n$. Set $\psi = (\psi_1,...,\psi_n)$. Consider the {\it generalized Bishop equation} 
\begin{eqnarray}
\label{Bishop1}
u(\zeta) = r(u(\zeta),T(u)(\zeta) + c) + t\psi(\zeta), \,\, \zeta \in b\D ,
\end{eqnarray}
where $c \in \R^m$ and $t = (t_1,...,t_m)$, $t_j \ge 0$, are real parameters. It follows by the implicit function theorem that this equation admits a unique solution 
$u(c,t) \in C^s(b\D)$ depending smoothly on the parameters $(c,t)$. Consider now the complex discs $f(c,t)(\zeta) = P(u(c,t)(\zeta) + iT(u(c,t))(\zeta))$, where $P$ denotes the Poisson operator of harmonic extension to $\D$:

\begin{eqnarray*}
P(u)(\zeta) = \int_{-\pi}^{\pi} K_P(\zeta,t)  u(e^{it})dt
\end{eqnarray*}
Here $K_P$ denotes the Poisson kernel
\begin{eqnarray*} 
K_P(\zeta,t) =  \frac{1}{2\pi} \frac{1- \vert \zeta \vert^2}{\vert e^{it} - \zeta \vert^2}.
\end{eqnarray*}
 The map $(c,t) \mapsto f(c,t)(0)$ (the centers of discs) is of class $C^s$.  Every disc is attached to $E$ along the upper semicircle. It is easy to see that this family of discs fills the wedge $W_\delta(\Omega,E)$ when $\delta > 0$ and a neighbourhood $\Omega$ of the origin are chosen small enough. Indeed, this is immediate when the function $r$ vanishes identically (i.e., $E = i\R^m$), while the general case follows by a 
small perturbation argument.

The detailed proof is contained in many works so I skip it (see, for example \cite{Su}).

\section{Proof of the main result}

In this section we prove Theorem \ref{MainTheo1}. For convenience of readers, we recall the general approach to the holomorphic mappings 
boundary value regularity problem,  and then explain how to modify this construction in order to  impove the regularity. In what follows we use the notation $C$, $C_1$,...
for positive constants. Their values may change from line to line.

\subsection{General construction}

The proof of Theorem \ref{MainTheo1} essentialy is based on estimates for the Kobayashi-Royden metric from \cite{ChSu}. However, they are quite different from the classical ones.



One of the first results on the boundary behaviour of holomorphic mappings (see \cite{Pi2}) is the following 

\begin{prop}
\label{ExtTheo1}
Let $f: \Omega_1 \to \Omega_2$ be a proper holomorphic mapping between two strictly pseudoconvex domains in $\C^n$ with boundaries of class $C^2$.
Then $f$ extends to $\overline\Omega_1$ as a H\"older 1/2-continuous mapping.
\end{prop}


The proof is based on the estimate of the Kobayashi-Royden metric given by Proposition \ref{KobProp1} and the Hopf Lemma.

\bigskip

In order imporove the boundary regularity of a mapping $f$, recall the construction due to S.Pinchuk and S.Khasanov. Let $\Omega$ be a strictly pseudoconvex domain of class $C^k$, $k \ge 2$. 

 If $\Omega = \rho < 0$, then explicitely we have 
$$H_p(b\Omega) = \{ v \in \C^n: \sum_{j = 1 }^n  \frac{\partial \rho}{\partial z_j}(p)v_j = 0 \}.$$
Similarly, denote by $H(b\Omega)$ the holomorphic tangent bundle of $b\Omega$, with the fiber $H_p(b\Omega)$ over a point $p \in b\Omega$. Every holomorphic tangent space is a complex hyperplane in $\C^n$ and can be viewed as a point of the complex projective space $\C\PP^{n-1}$. Therefore, the holomorphic tangent bundle is a real submanifold of dimension $2n-1$ and of class $C^{k-1}$ in $\C^n \times \C\PP^{n-1}$. This is well-known (and easy to check) that this manifold is totally real  when $b\Omega$ is strictly pseudoconvex
\cite{PiKh}.

Assume that $(\partial \rho / \partial z_n)(p) \neq 0$. Then 
$$H_p(b\Omega) = \{v \in \C^n:  v_n = w_1 v_1 + ...w_{n-1}v_{n-1}, w_j = \rho_{z_j}(p)/\rho_{z_n}(p) \}$$

Set $\phi_j(z) = \rho_{z_i}(z)/ \rho_{z_n}(z)$, $j= 1,...,n$. If $(w_1,...w_{n-1})$ are local coordinates in $\C\PP^{n-1}$  near $H_p(b\Omega)$, then in a neighborhood of 
$(p, H_p(b\Omega))$ the bundle $H_p(b\Omega)$ is defined by the equation 

\begin{eqnarray}
\label{bundle}
H(b\Omega) = \{ (z,w) \in \C^n \times \C^{n-1}: \rho(z) = 0, w_j = \phi_j(z) \}
\end{eqnarray}
and is the graph over $b\Omega$. Thus, the domain

\begin{eqnarray}
\label{Owedge}
W(\Omega) = \{ (z,w) \in \C^n \times \C^{n-1}: \rho(z) < 0 \}
\end{eqnarray}
is a domain with $C^k$-boundary in $\C^n \times \C^{n-1}$, containing in its boundary a $(2n-1)$- dimensional totally real submanifold $H(b\Omega)$ of regularity $C^{k-1}$.

\bigskip

The crucial step for the approach of Pinchuk-Khasanov is the following. Let $f: \Omega_1 \to \Omega_2$ be a proposer holomorphic mapping between two strictly pseudoconvex domains. We already know that $f$ is of class $C^{1/2}(b\Omega_1)$; in particular, $f$ is defined on the boundary $b\Omega$. We can define the lift $F$ of $f$ by to
$\C^n \times \C\PP^{n-1}$ by
$$F(z,P) = (f(z), df(z)P)$$
Here $P$ is viewed as a hyperplane in $\C^n$ and $df(z)P$ is its image under the tangent map $df(z)$ (recall that a proper holomorphic mapping between strictly pseudoconvex domain has a nonvanishing gradient). The map $F$ is defined on $W(\Omega_1)$ and takes it to $W(\Omega_2)$.

The following key result is due to Pinchuk-Khasanov.
\begin{lemma}
\label{PinLemma}
The map $F$ extends continuously on $W(\Omega_1) \cup H(b\Omega_1)$. The extension (again denoted by $F$) satisfies $F(z,H_z(b\Omega_1)) = (f(z),H_{f(z)}(b\Omega_2))$
for any $z \in b\Omega$. That is, $F(H(b\Omega_1)) \subset H(b\Omega_2)$.
\end{lemma}
The proof is based on the scaling method which in turn is based on the estimates of the Kobayashi-Royden metric  as in  Proposition \ref{KobProp1}. Further argument of Pinchuk-Khasanov is based on a smooth version of the edge-of-the-wedge theorem for asymptotically holomorphic functions. This argument requires the $C^s$ regularity with $s > 2$ for both 
$b\Omega_1$ and $b\Omega_2$. This is the reason why the $C^2$-case requires a different approach.

\subsection{Improving of regularity}

Consider a wedge-type domain $W$ given by (\ref{wedge}) with generic edge $E$ given by (\ref{edge}). Also, consider the shrinked wedge $W_\delta$ defined by (\ref{shrinked}).



We assume that the defining functions $\phi_j$ are of class $C^{1 + \alpha}$ with $\alpha > 0$,  and 
$\partial \phi_1 \wedge ... \wedge \partial \phi_m \neq 0$ in a neighborhood of $E$. %


Note that there exists a constant $C > 0$ such that for every point $z \in W_\delta$ one has 
\begin{eqnarray}
\label{dist}
C^{-1} dist(z, b W) \le dist(z, E) \le C dist (z,b W)
\end{eqnarray}

Now we can prove the following


\begin{thm}
\label{MainTheo2}
Let $N$ be an $n$-dimensional totally real manifold of class $C^1$ in  $\C^n$. Let also $W$ be a wedge in $\C^m$. Consider a holomorphic map   $f: W \to \C^n$ such that $f$ is continuos on $W \cup E$ and $f(E) \subset N$.  Then,  for every $\delta > 0$,  and every $\alpha < 1$, the map $f$ extends to a H\"older $\alpha$-continuous mapping on $W_{\delta} \cup E$. 
\end{thm}

In the case where $m=1$, that is $W$ is the unit disc, a much more general result is obtained in  \cite{ChSu}. Here we adapt the proof of \cite{ChSu} to our situation.

We begin the proof of Theorem with the following well-known lemma.

\begin{lemma}
\label{Poisson}
Let $\phi$ be a non-negative subharmonic function in $\D$ such that $\phi(\zeta) \to 0$ as $\zeta$ tends to an open 
arc  $\gamma \subset b\D$. Then, for every compact subset $K \subset \D \cup \gamma$, there exists a constant $C_K > 0$
such that $\phi(\zeta) < C_K(1 - \vert \zeta \vert)$ for any $\zeta \in K \cap \D$.
\end{lemma}
Let $V$ be a neighborhood of $\gamma \cap K$ such that $V’ = V \cap \D$ is simply connected and $\phi < 1$ in $V’$, and $g: \D \to V’$
be a conformal mapping.Then, by the reflection principle, $g^{-1}$ extends holomorphically across $\gamma$. Hence, replacing $\phi$ 
by $\phi \circ g$, we reduce the question to the case of a function that is uniformly bounded in $\D$. But then the assertion follows by an obvious estimate of the Poisson kernel.
Indeed, it suffices to consider the case where $\phi(e^{it})$ vanishes on the arc $\vert t \vert < \tau$ with some $\tau > 0$. For $\vert \arg \zeta \vert \le \tau/2$ and $\vert t \vert > \tau$ the Poisson kernel admits the estimate
$$K_P(\zeta,t) \le \frac{1}{\pi}\frac{1- \vert \zeta \vert}{\vert e^{it} - \zeta\vert^2} \le \frac{1}{\pi}\frac{1-\vert\zeta\vert}{\vert e^{i\tau} - e^{i\tau/2}\vert^2}.$$
But we have 
$$\phi(\zeta) \le \int_{-\pi}^{\pi} K_P(\zeta,t) \phi(e^{it})dt.$$
This gives the desired estimate and  proves lemma.

\bigskip

As a corollary we have 
\begin{lemma}
\label{Poisson2Lemma}
Let $\psi$ be a non-negative plurisubharmonic function on $W$ such that $\psi(z) \to 0$ as $z \to E$. Then for every fixed $\delta > 0$, there exists a constant $C = C_\delta > 0$
such that $\psi(z) \le C_\delta dist(z,E)$ for any $z \in W_\delta$.
\end{lemma}
For the proof consider the family $(h_t)$ of complex discs given by Lemma \ref{discs}. Then it suffices to apply Lemma \ref{Poisson}  to $\psi \circ h_t$  since all constants are uniform with respect to $t$.

\bigskip


By Lemma \ref{Poisson2Lemma}, there exists a constant $C_\delta > 0$ such that, for any $z \in W_\delta$,
one has 
\begin{eqnarray}
\label{Poisson2}
\rho \circ f (z) \le C_\delta dist(z,E)
\end{eqnarray}

We will need the result of Harvey-Wells \cite{HW}. Since $N$ is a totally real manifld of class $C^1$, there exists a non-negative function $\rho$ of class $C^2$ and strictly plurisubharmonic in a neеghborhhod of $N$, such that $N = \rho^{-1}(0)$. Furthermore, for each $\theta \in ]0,1[$, the function $\rho^{\theta}$ remains plurisubharmonic in a neighborhood of $N$.

Let $a \in E$ and $f(a) = p \in N$.  We can assume that $p = 0$ and there exists $\varepsilon > 0$ such that $\rho - \varepsilon \vert z \vert^2$ is plurisubharmonic on the ball $3 \B \subset D$.

\begin{lemma}
\label{LemKobEst1}
There exists a constant $A > 0$ with the following properry: if  $z \in W_{\delta}$ is a point such that $f(z) \in \B$, then
$$\parallel  df(z) \parallel \le A dist(z,E)^{-1/2}$$
(here we use the operator norm of the tangent map).
\end{lemma}
\proof Set $d = dist(z,E)$. By (\ref{dist}) there exists a constant $C > 0$ (independent of $z$) such that the ball $z + dC\B$ is contained in $W$.
 Then it follows from (\ref{Poisson2}) that the image $f(z+ dC\B)$ is contained in the domain $D_d =  \{ w \in D: \rho(w) < 2 C_1d \}$ if a constant $C_1 >0$ is suitably chosen. Note that  the strictly plurisubharmonic function $u_d(w): = \rho - 2 C_1d$ is negative on $D_d$. By Proposition \ref{ProLoc1}  there exists a constant $M > 0$, independent of $d$, such that, for any $w \in D \cap \B$ and any vector $\xi \in \C^n$, one has 
$$F_{D_d}(w,\xi) \ge M \vert \xi \vert \vert u_d(w) \vert^{-1/2}.$$
On the other hand, for the Kobayashi-Royden in  the ball $z + dC\B$, we have $F_{z+ d\B}(z,\tau) = \vert \tau \vert/dC$ for any vector $\tau \in \C^m$. Since the Kobayashi-Royden metric is decreasing under holomorphic mappings, one has
$$M \vert df(z)   \tau \vert \vert u_d(f(z)) \vert^{-1/2} \le F_{D_d}(f(z),df(z)\tau) \le F_{z + dC\B}(\zeta,\tau) = \vert \tau \vert/dC.$$
Therefore, $\vert df(z) \tau \vert \le M^{-1} \vert u_d(f(z)) \vert^{1/2}\vert \tau \vert/dC$. As $-2C_1 d \le u_d(f(\zeta)) < 0$, we have $\vert u_d(f(\zeta)) \vert^{1/2} \le (2C_1 d)^{1/2}$.
This implies the desried estimate and proves lemma.

\bigskip

Now it follows from Lemma \ref{LemKobEst1}  that  $f$ is H\"older $1/2$ continuous on $W_\delta \cup E$  via an integration argument which is a variation of the classical Hardy-Littlewood theorem  (cf. \cite{Ch}).


\bigskip



Let us improve the regularity. By Harvey-Wells theorem \cite{HW}, $\rho^\theta$ remains plurisubharmonic for every $1/2 < \theta < 1$. The composition $\rho^\theta \circ f$ is defined in a neighborhood $W \cup E$. 

Applying Lemma \ref{Poisson2Lemma} to this function, we obtain that $\rho \circ f (z) \le C dist(z,E))^{1/\theta}$ for each $z\in W \cup E$.

Now we simply repeat the former argument. Let $z \in W_\delta$ be sifficiently close to $a$ and $d = dist(z,E)$. The image $f(z + dC\B)$ is contained in the domain $D_d = \{ w \in D: u_d(w) = \rho(w) - 2 C_1 d^{1/\theta} < 0 \}$. Repeating the proof of Lemma \ref{LemKobEst1}, we obtain that $\parallel df(z) \parallel \le M^{-1}\vert u_d(f(\zeta)) \vert^{1/2}/dC$.
Since $-2C_1d^{1/\theta} \le u_d(f(\zeta)) < 0$, we conclude as above that that $\parallel df(\zeta) \parallel \le A dist(z,E)^{1/2\theta - 1}$ in a neighborhood of $a$ in $W$. Hence, 
$f$ is H\"older $1/2\theta$-continuous on $W_\delta \cup E$. This proves Theorem \ref{MainTheo2}.

\bigskip

Now we conclude the proof of Theorem \ref{MainTheo1}. It follows directly from Theorem \ref{MainTheo2}. It suffices to take $E = H(b\Omega_1)$ and $N = H(b\Omega_2)$ and to apply this theorem to the lift $F(z,P)$ of the map $f$. By Lemma \ref{PinLemma} the map $F$ satisfies the assumption of Theorem \ref{MainTheo2}. Since $H(b\Omega_1)$ is the graph over $b\Omega_1$, it is easy to choose a wedge $W \subset W(\Omega_1)$  (see (\ref{Owedge}))  with the edge $H(b\Omega_1)$, such that for any $\delta > 0$ the projection of the  shrinked wedge $W_\delta$ on $\C^n$ coincides with $\Omega$. Since $F$ is of class $C^{\alpha}(W_\delta \cup H(b\Omega_1))$, the initial map $f$ is of class $C^{\alpha}(\overline\Omega_1)$ for 
each $\alpha < 1$. 

\section{Further results}
A natural question arising is whether the conclusion of Theorem \ref{MainTheo1} still remains true when the boundary of $\Omega_1$ is only of class $C^2$. The answer is affirmative, 
but the proof is much more technical. We sketch here two possible approaches.

First of all we note that the assumption of $C^{2 + \varepsilon}$ regularity of $b\Omega_1$  allows to use Theorem \ref{MainTheo2} because  this condition guarantees that the regularity of $H(b\Omega_1)$  is $> 1$. The assumption that the regularity of $E$ is $> 1$ in Theorem \ref{MainTheo2}, in turn, is used in the proof only in order to establish Lemma \ref{Poisson2Lemma}. More precisely, this regularity assumption is necessary in order to apply Lemma \ref{discs} on gluing complex discs. Thus, it suffices to establish an analog of Lemma \ref{discs} for the case where the edge $E$ is exactly of class $C^1$.
This is possible but requires a much more careful analysis of the Bishop equation (\ref{Bishop1}).  Denote by $W^{k,p}(\D)$ the Sobolev classes of functions admitting  the Soblolev partial derivatives up to the order $k$  in  the class $L^p(\D)$. One can show that the equation (\ref{Bishop1}) admits a unique solution in $W^{1,p}(\D)$ for every $p > 2$. Furthermore, this solution is of class $C^1(\D)$ and depends $C^1$-smoothly on parameters in interior of  $\D$. This is sufficient in order to prove an analog of Theorem \ref{MainTheo2}. A detailed argument requires subtle results from the geometric measure theory.

The second approach uses the fact that in order to prove Theorem \ref{MainTheo1}, it suffices to establish Theorem \ref{MainTheo2} in a special case where the wedge $W$ coincides with $W(\Omega_1)$ and $E = H(b\Omega_1)$.  Technically it is more  appropriate to consider the holomorphic discs glued to the holomorphic tangent bundle  $H(b\Omega)$ (more precisely, to its projectivization) along the whole boundary (i.e. the whole unit circle $b\Omega$). This class of complex discs is well-known. This is exactly the stationary discs of Lempet \cite{Le} studied by several  authors. Suppose that $H(b\Omega_1)$ is defined by (\ref{bundle}). Then the holomorphic disc $z: \D \to \Omega_1$, $z: \zeta \mapsto z(\zeta)$ is called a stationary disc if it admits a holomorphic lift
$(z,w): \D \to \Omega_1 \times \C\PP^{n-1}$, $(z,w): \zeta \mapsto (z(\zeta),w(\zeta))$ with the boundary glued to $H(b\Omega_1)$ that is $(z,w)(b\D) \subset H(b\Omega_1)$. In the coordiantes this is equivalent to the Bishop type equation 
\begin{eqnarray*}
& &\rho(z(e^{i\theta})) = 0,\\
& &w_j  = \phi_j(z(e^{i\theta})), j= 1,...,n-1
\end{eqnarray*}
This is a non-linear Riemann-Hilbert type boundary value problem. This is well-known that the linearized problem is represented by a Fredholm operator which 
has positive partial indices and , therefore, is surjective. This  allows to solve the problem by the implicit function theorem (see, for example,  \cite{CoGaSu,SS}) in suitable function spaces such as $C^\alpha(\D)$ or $W^{1,p}(\D)$ with $p > 2$. In the case where 
$\Omega_1 = \B^n$ (which corresponds to the linearized Riemann-Hilbert problem), the stationary  discs belong to complex lines.
Consider a point $a \in b \Omega_1$. One can assume that $\B^n$ is tangent to $\Omega_1$ at $a$ up to the order 2. Let $b \in \Omega_1$ be a  point on the real inward normal to $b\Omega_1$ at the point $a$. 
 Consider  vector $v\in \C^n$ parallel to $H_a(b\Omega_1)$. 
 The unique stationary disc through $b$ in the direction $v$ is a small perturbation of a complex line  $L$ through $b$ in the direction $v$. Note that its lift glued to $H(b\Omega_1)$ is a large complex disc (the totally real manifold $H(b\Omega_1)$ does not admit small attached complex discs). These discs form a foliation which is a small deformation of a foliation of $\Omega_1$ near a boundary point $p \in L \cap b\Omega_1$ by a family of complex lines parallel to  holomorphic tangent space of $b\Omega$. Every disc is of class $C^{\alpha}$ with $\alpha < 1$ and then the proof of Theorem \ref{MainTheo2} goes through.

Detailed presentation of these approaches will appear in a forthcoming paper.

\end{document}